\documentclass[journal]{IEEEtran}

\usepackage{amsmath}
\usepackage{amssymb}
\usepackage{amsfonts}

\usepackage{amsthm}
\usepackage[mathscr]{euscript}
\usepackage{bbm}
\usepackage{mathtools}
\usepackage{breqn}
\usepackage{setspace}
\usepackage{gensymb}
\usepackage{url}
\usepackage{lipsum}
\usepackage{soul,color}
\usepackage{graphics}
\usepackage{graphicx}
\usepackage{epstopdf}
\usepackage{subcaption}
\usepackage{siunitx}
\usepackage{enumitem}
\usepackage{algorithmicx}
\usepackage{algpseudocode}
\usepackage[Algorithm,ruled]{algorithm}
\usepackage{xcolor}
\usepackage{soul}
\usepackage{cite}
\usepackage{kantlipsum}

\DeclareMathOperator*{\minimize}{minimize}

\soulregister\cite7
\soulregister\ref7
\soulregister\pageref7
\soulregister\eqref7

\pdfoutput=1

\begin{document}
	
	\title{
	Energy Consumption Optimization for Autonomous Vehicles via Positive Control Input Minimization}
	
	\author{Andreas Hadjigeorgiou and Stelios Timotheou

		\thanks{This work is supported by the European Union (i. ERC, URANUS, No. 101088124 and, ii. Horizon 2020 Teaming, KIOS CoE, No. 739551), and the Government of the Republic of Cyprus through the Deputy Ministry of Research, Innovation, and Digital Strategy. Views and opinions expressed are however those of the author(s) only and do not necessarily reflect those of the European Union or the European Research Council Executive Agency. Neither the European Union nor the granting authority can be held responsible for them.}
		\thanks{The authors are with the Department of Electrical and Computer Engineering and the KIOS Research and Innovation Center of Excellence, University of Cyprus {\tt\small \{hadjigeorgiou.g.andreas, timotheou.stelios\}@ucy.ac.cy}.}
		
	}

	\maketitle

	\begin{abstract}

Autonomous vehicles (AVs) present a unique opportunity to improve the sustainability of transportation systems by adopting eco-driving strategies that reduce energy consumption and emissions. This paper introduces a novel surrogate model for energy and fuel consumption that minimizes Positive Control Input (PCI). Unlike conventional objectives such as squared acceleration, which often misrepresent actual energy usage, PCI provides a more accurate and optimization-friendly alternative. Building on PCI, we propose ECO+, a convex, time-based trajectory optimization framework that ensures safety and passenger comfort while optimizing energy use for AVs approaching an intersection. To improve computational efficiency, quadratic resistive forces are approximated using piecewise affine segments, resulting in a linear programming formulation. ECO+ is validated using empirical fuel and electric energy models and benchmarked against established optimization strategies, including a state-of-the-art nonlinear solver. Simulation results demonstrate that ECO+ consistently outperforms baseline methods in reducing energy consumption, even under strict comfort constraints and in scenarios involving a leading vehicle. Moreover, initializing a nonlinear solver with ECO+ yields only marginal gains, indicating that ECO+ is effective as a standalone eco-driving strategy. These findings highlight ECO+ as a practical, scalable, and computationally efficient solution for enhancing the sustainability of autonomous urban mobility systems.

	\end{abstract}
	

	\IEEEpeerreviewmaketitle
	
		\section{Introduction}
        \label{sec:Introduction}

The urgent need for greener transportation arises from the necessity to mitigate the environmental impacts of conventional transportation systems, which are major contributors to greenhouse gas emissions (GHG), thereby accelerating climate change and degrading urban air quality. The transportation sector currently accounts for approximately 25\% of Europe’s total GHG emissions, with road transportation alone contributing around 70\% of these emissions. Cars, vans, trucks, and buses, which predominantly rely on fossil fuels, are responsible for the majority of this impact, consuming approximately 340 billion liters of petroleum fuel annually in Europe \cite{transportation_eceee}. Transitioning to sustainable transportation practices, such as the adoption of fuel-efficient technologies, alternative energy sources, and intelligent traffic management systems, is imperative to reduce emissions and enhance energy efficiency \cite{shah2021green}.

Despite advancements in vehicle efficiency, the increasing global population and expanding economies have led to a continuous rise in vehicle-miles traveled rather than a reduction. To counteract this trend, sustainable transportation solutions are being actively explored. One promising approach is \textit{eco-driving}, which seeks to minimize fuel consumption through energy-aware driving behaviors and optimized vehicle control.

Efforts to promote energy-aware driving behaviors rely on heuristic guidelines aimed at human drivers. These include recommendations such as shifting to higher gears early, anticipating traffic flow to minimize unnecessary braking, maintaining an energy-optimal speed, and turning off the engine during stops in congestion or at traffic lights \cite{D3}. Such approaches involve eco-driving training programs \cite{greek_bus} and software-based assistance systems that aid drivers in making energy-efficient decisions. These systems support drivers in eco-routing, helping them navigate through energy-efficient routes \cite{eco-routing}, and provide real-time feedback to encourage energy-efficient driving behaviors \cite{freeway_dynamic_advice, pEDAS}. Additionally, some approaches incentivize eco-driving through post-trip feedback, displaying comparative energy-related statistics to encourage long-term behavioral change \cite{D4}. Collectively, these strategies have demonstrated fuel savings of approximately 10--20\%.

While these approaches are effective, their success depends on driver compliance and behavioral adaptation, which can vary significantly among individuals \cite{fuel_survey}. Autonomous vehicles (AVs)  have emerged as a promising technology, offering a transformative solution by eliminating human variability and executing optimized energy-efficient driving strategies consistently.

Within this context, AVs can enhance eco-driving strategies, reducing fuel and energy consumption, and improving overall traffic efficiency \cite{CAV_better}. By leveraging advanced control algorithms and real-time data processing, AVs can play a pivotal role in optimizing transportation systems while aligning with global sustainability goals.

Extensive research has been conducted on mathematically rigorous eco-driving frameworks that account for various factors, including powertrain type (e.g., internal combustion engines \cite{C7}, hybrid systems \cite{D}, and electric vehicles \cite{C13}), road geometry (e.g., curvature and inclination \cite{kamal}), vehicle type (e.g., passenger cars \cite{C7} and heavy-duty vehicles \cite{C8,C9}) and alternative objectives (e.g. minimizing emmisions). Additionally, these frameworks consider different problem settings, such as single-vehicle eco-driving \cite{oneVehicle}, cooperative adaptive cruise control \cite{C}, and platooning \cite{platoonVehicles}, as well as key assumptions affecting vehicle dynamics, such as the influence of aerodynamic drag \cite{D}.

Formulating eco-driving problems within a mathematically rigorous framework often results in complex, non-convex optimization problems that pose significant computational challenges. Consequently, various solution methodologies have been proposed in the literature, including dynamic programming \cite{A12}, optimal control \cite{A13}, mixed-integer programming \cite{eco-driving_MIPNL}, and machine learning-based approaches \cite{AI_energy,A19}. However, each of these methods has inherent limitations. Dynamic programming and mixed-integer programming often suffer from the curse of dimensionality, making them impractical for real-time execution. Optimal control techniques require the derivation of analytical expressions, which is cumbersome and typically applicable only to simplified problem formulations. Machine learning approaches, while promising, necessitate extensive offline training with large datasets and do not provide formal guarantees regarding the quality of the generated solutions.

Acknowledging these limitations, researchers have explored heuristic strategies as an alternative. A particularly effective and popular approach is pulse-and-glide (PnG), a driving strategy characterized by alternating periods of acceleration and coasting \cite{lee2009vehicle}, which has been extensively studied in various contexts.  These include applications in mechanical step-gear AVs \cite{PnG_gear}, free-driving scenarios \cite{PnG_free}, car-following situations \cite{PnG_following}, mixed traffic conditions \cite{PnG_mixed}, and investigations into ride comfort trade-offs \cite{PnG_tradeoff}. Across these studies, PnG has demonstrated in simulations and experiments significant fuel savings, often exceeding 30\%. However, most of this research has primarily focused on highway scenarios, with limited attention given to urban environments.

Urban traffic conditions, particularly intersections, present considerable opportunities for fuel consumption reduction \cite{RiosSurvey}. Intersections serve as critical points where multiple traffic flows converge, and optimizing vehicle movements in these areas can yield substantial energy savings. This has led to a growing body of research on autonomous intersection control \cite{Latest_survey}, exploring various optimization objectives. However, relatively few studies explicitly consider eco-driving. Some research efforts have sought to minimize transient engine operation \cite{Malikopoulos2018,xu2022general,ahadjistimoITV}, though alternative optimization objectives may yield even greater fuel savings \cite{fuel_convex_concave}. Another approach involves minimizing the positive tractive force applied by the AV \cite{simos_stelios_CAV}; however, such formulations often rely on simplifying assumptions regarding feasibility and approximation accuracy.

Other studies \cite{Kamal_smooth,hult2020optimisation} assess the performance of intersection coordination controllers without directly incorporating fuel consumption optimization as an explicit control objective. A widely adopted strategy in eco-driving research is model predictive control (MPC) \cite{ahadjistimoITV,fuel_convex_concave,simos_stelios_CAV,Kamal_smooth,hult2020optimisation,C16,C40}, which is well-suited for real-time applications and handling large-scale nonlinear problems. Nevertheless, the optimization problems formulated within this framework are frequently non-convex, making them computationally challenging to solve. The standard practice in such cases is to employ off-the-shelf nonlinear solvers, which may result in suboptimal or poor-quality solutions, and whose high computational costs limit their applicability in real-time scenarios.

This work addresses the challenge of achieving real-time, scalable eco-driving performance for AVs approaching intersections, where computational efficiency is critical. Specifically, we consider the problem of minimizing the energy consumption of an AV as it traverses a prescribed road segment, subject to fixed initial conditions and terminal constraints on both time and speed. This scenario is frequently encountered in autonomous intersection management, where precise coordination and operational efficiency are essential \cite{surv_AIM}.

To tackle this problem, we employ a discrete-time MPC framework. Preliminary work by the authors focused exclusively on minimizing fuel consumption and introduced a customized difference-of-convex (DC) solution approach \cite{fuel_convex_concave}. In this paper, we extend both the formulation and the solution methodology. First, we incorporate additional considerations, including passenger comfort and interactions with a leading vehicle. Moreover, rather than directly optimizing the fuel or energy consumption models, we introduce a novel convex surrogate objective to enhance computational tractability. The proposed method is evaluated using two distinct models: one for fuel consumption and one for energy consumption.

The key contributions of this paper are as follows: 
\begin{itemize} 
\item Introduction of PCI (Positive Control Input), a convex and efficient surrogate objective function designed to capture energy/fuel consumption by minimizing only the positive control effort. The PCI surrogate function generalizes across different vehicle types without requiring calibration. 
\item Development of ECO+, a convex optimization framework incorporating the PCI surrogate consumption model, vehicle dynamics, resistive forces, comfort constraints, and safety requirements for AVs approaching intersections. A piecewise affine approximation of quadratic constraints is proposed to transform the problem into a linear program suitable for real-time implementation. 
\end{itemize}

ECO+ is benchmarked against several established optimization strategies, including a state-of-the-art nonlinear solver. Extensive simulations across a range of traffic scenarios demonstrate that it consistently outperforms alternative methods in reducing both fuel and energy consumption.

Overall, the proposed solution approach delivers real-time, scalable, and high-quality solutions to the eco-driving problem by leveraging the novel PCI model and a linear-programming-based MPC framework,  enhancing the efficiency and sustainability of autonomous vehicle operations.

The remainder of the paper is structured as follows. Section \ref{sec:formulation} introduces the problem setup and formulates the optimization problem. Section \ref{sec:EcoModels} details the fuel and energy consumption models used to evaluate eco-driving performance. Section \ref{sec:Methodology} presents the proposed solution approach. Section \ref{Sec:Results} introduces the benchmark methods and presents a comprehensive performance evaluation across various scenarios. Finally, Section \ref{sec:conclusion} summarizes our key findings and conclusions.

\section{Problem Statement}
\label{sec:formulation}

\begin{figure}[t]
	\centering
	\includegraphics[width=\columnwidth]{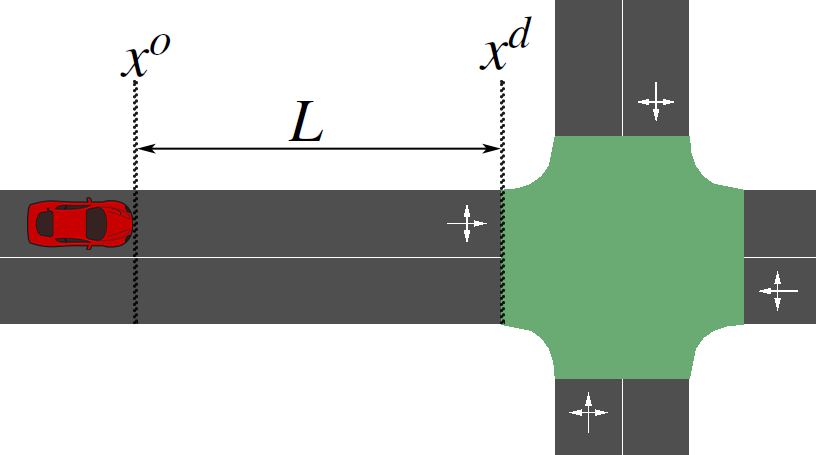}
	\caption{Illustration of an AV approaching an intersection.}
	\label{fig:Problem_illustration}
\end{figure}

This work considers an AV traveling from an origin, $x^o$, to a specified destination, $x^d$. The objective of the AV is to minimize fuel/energy consumption while covering a total distance of $L$ meters. The AV starts its journey at time $t^o$ with an initial speed of $v^o$ and aims to reach the destination at time $t^d$ with a final speed of $v^d$, resulting in a total travel time of $t^m = t^d - t^o$. As illustrated in Figure~\ref{fig:Problem_illustration}, this scenario represents a typical case where an AV must pass through an unsignalized intersection at a specified time and speed as instructed by an intersection controller. This intersection controller ensures the safe and seamless passage of multiple AVs through the intersection without requiring them to stop, as described in \cite{ahadjistimoITV}.

The travel time and acceleration profile, as explored in \cite{hadjigeorgiou2024optimizing}, play a critical role in determining consumption, particularly given the cumulative effect of numerous vehicles traveling through an intersection. Thus, developing an optimization framework and an efficient methodology both in terms of computational cost and solution quality is highly desirable for achieving real-time, sustainable and coordinated traffic flow.

In this section, we state and formulate the optimization problem that serves as the foundation for achieving the objective of minimizing consumption. The proposed framework incorporates key factors such as vehicle dynamics, speed limits, passenger comfort, friction, and safety considerations. 

\subsection{AV Modeling}
To begin, the vehicle dynamics are modeled as a third-order system to accurately capture the relationship between control inputs and vehicle behavior, ensuring a realistic representation. Specifically,
\begin{subequations}\label{eq:dynamics}
\begin{align}
\dot{x}(t)=v(t),\quad \forall t\in [t^o,t^d], \label{eq:dynamicsA}\\
	 \dot{v}(t)=a(t),\quad \forall t\in [t^o,t^d], \label{eq:dynamicsB}\\
	 a(t)=u(t)-a^r(t),\quad \forall t\in [t^o,t^d], \label{eq:dynamicsC}\\
	 \dot{a}(t)=J(t),\quad \forall t\in [t^o,t^d] \label{eq:dynamicsD}
\end{align}
\end{subequations}
where $x(t)$, $v(t)$, $u(t)$, $a(t)$, $J(t)$ are the position, velocity, control input, acceleration, jerk at time $t$, respectively. The term $-a_r(t)$ represents the deceleration experienced by the vehicle due to opposing forces such as rolling resistance and air drag, as well as the acceleration or deceleration caused by the road slope. For small road segments, such as the one considered in this scenario, the road slope, $\theta$, can be considered constant. To include the boundary conditions  as discussed above, the following constraints must be enforced:
\begin{subequations}\label{eq:dynamics_boundary}
\begin{align}
	&\enspace  x(t^o)=0 , \quad x(t^d)=L,\\
	&\enspace  v(t^o)=v^o, \quad v(t^d)=v^d, \\
	&\enspace  u(t^o)=0,\quad u(t^d)=0.
\end{align}
\end{subequations}

While most formulations in the literature model vehicle dynamics as a second-order system, it is equally important for such formulations to account for passenger comfort. Incorporating jerk, the derivative of acceleration, into the model allows for the implementation of thresholds that ensure comfort. Specifically, this can be achieved by placing limits on the jerk, thereby keeping it within acceptable tolerance levels as 
\begin{equation}\label{eq:jerk_constraint}
\begin{aligned}
J^\min\leq J(t)\leq J^\max,\quad \forall t\in [t^o,t^d],
\end{aligned}
\end{equation}
where $J^\min$ and $J^\max$ are the lower and upper limits of jerk respectively. If comfort is not of concern then \eqref{eq:dynamicsD} and \eqref{eq:jerk_constraint} can be omitted.
Resistive forces are explicitly incorporated and their significance is analyzed in Section \ref{Sec:Results}. These forces include rolling resistance, aerodynamic drag, and the influence of road slope.

The rolling resistance force, \( F_R \), is typically modeled as \( F_R = \mu m g \cos(\theta) \), where \( \mu \) is the rolling resistance coefficient, \( m \) is the mass of the vehicle and \( g \) is the gravitational acceleration. Alternatively, more recent standards model \( F_R \) as a function of velocity \cite{SAE_J2452_201707}, expressed as \( F_R = c_0 + c_1 v \), where \( c_0 \) and \( c_1 \) are empirically derived constants.

The aerodynamic drag force is given by \( F_D = \rho C_D A \tilde{v}^2/2 \), where \( \rho \) is air density, \( C_D \) is the drag coefficient, \( A \) is the frontal area of the vehicle, and \( \tilde{v} \) is the velocity relative to air. For simplicity, we assume still air, so \( \tilde{v} = v \).

The gravitational force due to the slope of the road is \( F_W = m g \sin(\theta) \), which can either assist or resist motion depending on whether the slope is downhill or uphill, respectively.

These resistive forces can be consolidated into an equivalent deceleration term,
\begin{equation}\label{eq:resistive_and_slope_equation}
\begin{aligned}
a^r(t) = d_1 + d_2 v(t) + d_3 v^2(t),
\end{aligned}
\end{equation}
where \( d_1 \), \( d_2 \), and \( d_3 \) aggregate the effects of rolling resistance, aerodynamic drag, and slope using Newton’s second law.

Next, the velocity and acceleration must remain within physical limits or externally imposed constraints at all times, which may stem from infrastructure requirements, safety regulations, passenger comfort, or the vehicle's design specifications. To enforce this, it is imposed that
\begin{equation}
	\label{eq:speed_limit}
	\begin{aligned}
	& 0\leq v(t)\leq v^{\max}, \quad  \forall t\in [t^o,t^d],  \\
	&u_i^{\min}\leq u(t)\leq u^{\max},  \quad \forall t\in [t^o,t^d].
	\end{aligned}
\end{equation}

\subsection{Safety constraints}

The final constraint considered ensures the safety of the ego AV with respect to a preceding vehicle traveling along a known trajectory, \( x^f(t) \). To maintain a safe following distance, it is imposed that
\begin{equation}
\label{eq:safety_distance} 
x^f(t) - x(t) \geq \max\left(\delta, \left(v(t) - v^f(t)\right)t^g\right)
\end{equation}
where \( \delta \) denotes the minimum allowable distance between the ego and the preceding vehicle, \( v^f(t) \) is the velocity of the preceding vehicle, and \( t^g \) is the desired time gap. This constraint integrates two safety components. 

The first term, \( x^f(t) - x(t) \geq \delta \), ensures a fixed minimum distance to prevent collisions, particularly in scenarios involving sudden deceleration of the preceding vehicle. Although effective at low speeds, a fixed buffer may be insufficient at higher velocities \cite{time_Gap_ci}. Thus, the second term, \( x^f(t) - x(t) \geq (v(t) - v^f(t))t^g \), enforces a dynamic, speed-dependent time gap. This accounts for factors such as communication delays and control or sensing inaccuracies that become more critical at higher speeds \cite{time_Gap_ci}.

Together, these components provide comprehensive safety measures: the fixed distance protects against abrupt stops at lower speeds, while the time-gap formulation enhances safety at higher speeds through dynamic spacing based on relative velocity.

\subsection{Objective}

To optimize energy consumption without directly minimizing the consumption models, which are typically non-convex and computationally expensive, we propose the PCI surrogate function defined as
\begin{equation}
\label{eq:heuristic_pos_acc}
\int_{t^o}^{t^d} u^{+}(t),
\end{equation}
where $u^{+}(t) = \max(u(t), 0)$ is the positive control input.

The PCI surrogate consumption model focuses exclusively on minimizing positive control effort, which directly corresponds to the energy-intensive phases of vehicle operation. By targeting these high-energy segments, PCI effectively reduces overall energy expenditure by inducing sparsity on the control input.

The model is inspired by pulse-and-glide strategies widely discussed in the literature, as well as extensive simulation-based experimentation conducted by the authors. Through iterative testing and analysis, the formulation was refined to identify acceleration patterns that consistently yield lower energy consumption, ultimately motivating the decision to minimize positive control input directly.

In addition to its simplicity, the PCI model is convex, leading to optimization problems that are computationally efficient and solvable in real time, and broadly applicable across vehicle types. While it does not require calibration, its performance assumes that resistive forces are modeled with reasonable accuracy, which are typically parameterized in vehicle dynamics and can be estimated with standard techniques.
As such, it is a practical and scalable objective function that can be seamlessly integrated into existing eco-driving frameworks to enhance energy efficiency.

\subsection{Problem Formulation}

In sum, the problem of determining the control input, $u(t)$,  that minimizes consumption can be expressed as
\begin{subequations} 	\label{eq:general_problem_continuous}
	\begin{align}
	\displaystyle{\minimize_{u(t)}}&\enspace \int_{t^o}^{t^d} u^{+}(t) \label{eq:general_problem_continuousA}\\
	\text{subject to:}&\enspace \eqref{eq:dynamics}-\eqref{eq:safety_distance},
	\end{align}
\end{subequations}
where $u(t)$ is the main optimization variable and $x(t)$, $v(t)$, $a(t)$, $J(t)$ are state variables.

This formulation is inherently non-convex due to the resistive force constraint \eqref{eq:resistive_and_slope_equation}. The solution approach 
is presented in detail in Section~\ref{sec:Methodology}.

\section{Consumption models}
\label{sec:EcoModels}

To effectively evaluate eco-driving techniques, a suitable energy or fuel consumption model capable of predicting instantaneous consumption on a second-by-second basis is essential.

The survey in \cite{fuel_survey} classifies fuel consumption models into three categories: white-box, black-box, and grey-box models. White-box models, such as \cite{SAERENS20091582}, rely on a detailed understanding of the system being modeled. However, they often require parameters that are difficult or impossible to obtain, making them impractical for many researchers. In contrast, black-box models, such as \cite{Malikopoulos2013}, are purely data-driven, relying solely on input-output relationships without incorporating physical principles. While these models can be effective, they lack interpretability and generalizability. Grey-box models serve as an intermediate approach, combining both system insights and experimental data. By integrating physics-based principles with empirical observations, grey-box models strike a balance between accuracy and complexity, making them a practical choice for eco-driving applications \cite{fuel_survey}.

Considering these factors, we adopt two widely used grey-box models: one tailored for energy consumption in battery electric vehicles and another for fuel consumption in internal combustion engine vehicles. Together, these models provide a reliable framework for analyzing eco-driving strategies across both vehicle types.

 \subsection{CPEM Model}
\label{subsec:CPEM}

\begin{figure}[t!]
	\centering
	\begin{subfigure}[b]{\columnwidth}
		\centering
		\includegraphics[width=1\linewidth]{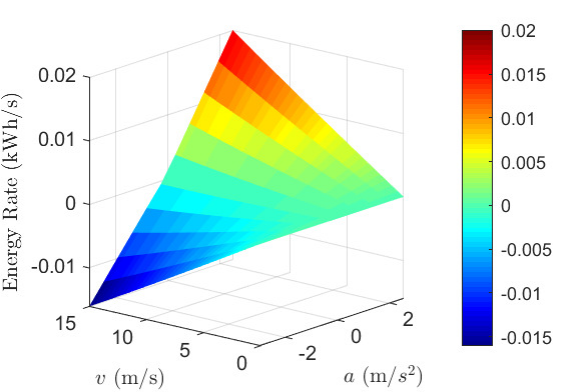}
		\caption{CPEM model}
		\label{fig:CPEM_model}
	\end{subfigure}%
	
	\begin{subfigure}[b]{\columnwidth}
		\centering
		\includegraphics[width=1\linewidth]{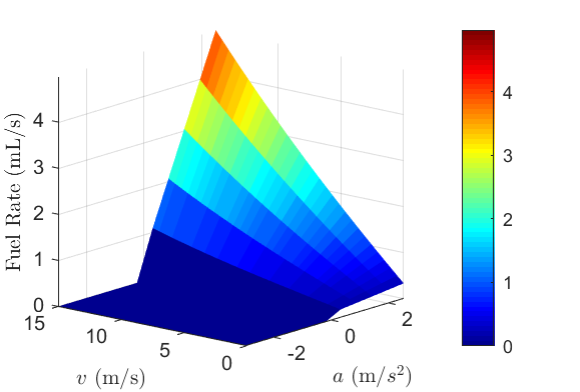}
				\caption{KMMK model}
		\label{fig:Kamal_model}
	\end{subfigure}
	\caption{Fuel and energy consumption rates for varying acceleration and velocity for each model.}
	\label{fig:Model_plot}
\end{figure}

One of the most widely used energy consumption models for electric vehicles is the Comprehensive Power based Electric vehicle energy consumption Model (CPEM), introduced in \cite{CPEM}. This model approximates energy consumption (EC) as a function of acceleration and velocity. The authors first consider the power at the wheels, $P_{\text{w}}(t)$, which is calculated using the following formula
\begin{equation}
\label{eq:CPEM_p1}
\begin{aligned}
P_{\text{w}}(t) = \Big(m a(t) + m g \cos{\theta} \frac{C_r}{1000} (c_1 v(t) + c_2) \\
+ \frac{1}{2} \rho A C_D v^2(t) + m g \sin{\theta} \Big) v(t),
\end{aligned}
\end{equation}
where $c_1$, $c_2$, and $C_r$ represent rolling resistance parameters that vary based on the road surface type, road conditions, and vehicle tire type. The formula is derived by accounting for the main forces acting on the vehicle. Using the power at the wheels, the power at the electric motor is given by
\begin{equation}
\label{eq:CPEM_p2}
P_{\text{em}}(t) = \frac{P_{\text{w}}(t)}{\eta_{\text{d}} \cdot \eta_{\text{em}}},
\end{equation}
where $\eta_{\text{d}}$ is the driveline efficiency and $\eta_{\text{em}}$ is the efficiency of the electric motor. Based on $P_{\text{em}}(t)$, the instantaneous power required or generated by the electric motor is expressed as
\begin{equation}
\label{eq:CPEM_p3}
P(t) =
\begin{cases}
P_{\text{em}}(t), & P_{\text{w}}(t) \geq 0,\\
P_{\text{em}}(t) \cdot \eta_{\text{rb}}, & P_{\text{w}}(t) < 0,
\end{cases}
\end{equation}
where $\eta_{\text{rb}}$ is the regenerative braking efficiency. 
The total energy consumption can then be calculated by integrating the required power, $P(t)$, over the travel time as
\begin{equation}
\label{eq:CPEM_p4}
\text{EC} = \frac{1}{3600000} \int^{t^m}_0 P(t) \cdot \eta_{\text{b}} \, dt,
\end{equation}
where EC represents the total energy consumption in kWh, and $\eta_{\text{b}}$ is the battery efficiency. The energy consumption rate, $P(t) \cdot \eta_{\text{b}}$, is illustrated in Figure \ref{fig:CPEM_model}. The figure shows that the energy consumption rate increases with increasing acceleration and velocity. Conversely, during deceleration, energy is regenerated, although at a reduced rate due to the limited efficiency of regenerative braking.

 \subsection{KMMK Model}
\label{subsec:KMMK}

For fuel consumption, we adopt the widely used model developed by Kamal, Mukai, Murata, and Kawabe (KMMK) in \cite{kamal}. KMMK model is comprised of a polynomial function of acceleration and speed to approximate instantaneous fuel consumption by considering the forces acting on a vehicle and the engine efficiency characteristics. First, the fuel consumption is approximated for a cruising vehicle (without acceleration or braking) using a curve-fitting process, resulting in
\begin{equation}
\label{eq:Fuel_kamal_1}
\begin{aligned}
f_{\text{cruise}} = c_0 + c_1 v(t) + c_2 v^2(t) + c_3 v^3(t),
\end{aligned}
\end{equation}
where $c_0$ to $c_3$ are curve-fitting coefficients. Using a similar procedure, the fuel consumption due to positive acceleration can be calculated as
\begin{equation}
\label{eq:Fuel_kamal_2}
\begin{aligned}
f_{\text{accel}} = \hat{a}(t) \left( c_4 + c_5 v(t) + c_6 v^2(t) \right),
\end{aligned}
\end{equation}
where $c_4$ to $c_6$ are curve-fitting coefficients and $\hat{a}(t)=a_v(t) +a_{\theta}(t)$, where $a_v(t)=-C_D \rho Av^2(t) /(2M)-\mu g\cos{(\theta)}-g\sin{(\theta)}+u(t)$ is the apparent acceleration and $a_\theta=g\sin{\theta}$ is the acceleration internally required to counteract the deceleration due to the road slope. 

The total fuel consumption is obtained by integrating the fuel consumption rates for cruising and accelerating over the duration of the trip. The total fuel consumption for the KMMK model is given by
\begin{equation}
\label{eq:Fuel_kamal}
\begin{aligned}
FC = \int_{t^o}^{t^d} \left( f_{\text{cruise}} + f_{\text{accel}} \right) \mathbbm{1}_{(0, u^{\max}]}(u(t)) \, \mathrm{d}t,
\end{aligned}
\end{equation}
where $FC$ is measured in milliliters (mL)
. The indicator function, $\mathbbm{1}_{(0, u^{\max}]}(u(t))$\footnote{$\mathbbm{1}_{ (a, b] }(u(t))$ is the indicator function which is equal to 1 if $ u(t)\in(a, b] $ and 0 otherwise.}, ensures that fuel consumption is zero when the control input is less than or equal to zero.

The fuel consumption rate function of the KMMK model is illustrated in Figure \ref{fig:Kamal_model}. Similar to the CPEM model, fuel consumption increases with higher acceleration and velocity. However, in contrast to CPEM model, where power is stored in the battery when the control input is less than or equal to zero, the KMMK model assumes no fuel consumption in such scenarios.

\section{Solution Approach}
\label{sec:Methodology}

In this section, we present the solution approach to the optimization problem \eqref{eq:general_problem_continuous}. To balance computational efficiency and solution quality, the continuous-time problem is first discretized. The resulting formulation is then transformed into an equivalent convex form, followed by linearization of the quadratic constraints.

\subsection{Discretization}
\label{subsec:Discretization}

Euler's forward method provides a first-order approximation of the continuous-time dynamic system $\dot{y} = f(t, y)$ via the discrete-time update rule
\begin{equation}\label{eq:dynamics_discrete}
y_{i+1} = y_i + f(t_i, y_i) \Delta T,
\end{equation}
where $\Delta T$ denotes the time discretization step. Applying this scheme to problem~\eqref{eq:general_problem_continuous} yields the following discrete-time optimization problem
\begin{subequations}
\label{eq:general_problem_discrete}
\begin{align}
\displaystyle{\minimize_{\mathbf{u}}}&\enspace \sum_{i\in\mathcal{I}} u^{+}_i  \Delta T\ \label{eq:general_problem_discreteA} \\
\text{subject to:}&\enspace x_{i+1}=x_i+ \Delta T\  v_i \quad \forall i\in\mathcal{I},\label{eq:general_problem_discreteB}\\
&\enspace v_{i+1}=v_i+ \Delta T\  a_i \quad \forall i\in\mathcal{I}, \label{eq:general_problem_discreteC}\\
&\enspace a_{i}=u_i- a^r_i \quad \forall i\in\mathcal{I}, \label{eq:general_problem_discreteD}\\
&\enspace a^r_{i}= d_1 + d_2 v_i + d_3 v_i^2 \quad \forall i\in\mathcal{I}, \label{eq:general_problem_discreteE}\\
&\enspace J_{i}\Delta T=a_{i+1}-  a_i \quad \forall i\in\mathcal{I}, \label{eq:general_problem_discreteF}\\
&\enspace 0\leq v_i\leq v_i^{\max} \quad \forall i\in\mathcal{I^+},\label{eq:general_problem_discreteG}\\
&\enspace u^{\min}\leq u_i\leq u^{\max} \quad \forall i\in\mathcal{I^+}, \label{eq:general_problem_discreteH}\\
&\enspace J^{\min}\leq J_i\leq J^{\max} \quad \forall i\in\mathcal{I^+}, \label{eq:general_problem_discreteI}\\
&\enspace x_0=0 , \quad x_H=L,	\label{eq:general_problem_discreteJ} \\
&\enspace  v_0=v^o, \quad v_H=v^d,	\label{eq:general_problem_discreteK}\\
&\enspace u_0=0,\quad u_H=0, \label{eq:general_problem_discreteL}\\
&\enspace x^f_i-x_i \geq \delta , \label{eq:general_problem_discreteM}\\
&\enspace x^f_i-x_i \geq  (u_i-u^f_i) t^g, \label{eq:general_problem_discreteN}
	\end{align}
\end{subequations}
where   $H= \left \lceil{t^m/\Delta T}\right \rceil$ with $\lceil\cdot\rceil$ being the ceiling function,  $\mathcal{I}=\{ 0,1,\dots,H-1 \}$ and $\mathcal{I^+}=\{ 0,1,\dots,H \}$. Note that the main optimization variables is $\mathbf{u}=\{u_0,u_1,\dots,u_{H} \}$ and the state ones are
$\mathbf{x}=\{x_0,x_1,\dots,x_{H} \}$, $\mathbf{v}=\{v_0, v_1,\dots,v_{H} \}$, $\mathbf{a}=\{a_0,a_1,\dots,a_{H} \}$ and $\mathbf{J}=\{J_0,J_1,\dots,J_{H} \}$.

In the following subsection, we detail the transformation of the problem into an equivalent convex form and the reformulation of the quadratic constraint~\eqref{eq:general_problem_discreteE} as a set of linear inequalities. This enables efficient numerical solution of the optimization problem~\eqref{eq:general_problem_discrete}.

\subsection{Convex Reformulation}

Constraint~\eqref{eq:general_problem_discreteE} is a quadratic equality, which is non-convex and cannot be handled directly by standard solvers. To address this, we first combine constraints~\eqref{eq:general_problem_discreteD} and~\eqref{eq:general_problem_discreteE} into a single expression, resulting in
\begin{equation}
\label{eq:combined}
a_i + d_1 + d_2 v_i + d_3 v_i^2 = u_i, \quad \forall i \in \mathcal{I^+}.
\end{equation}
The optimization problem can now be written as
\begin{equation}
\begin{aligned}
\displaystyle{\minimize_{\mathbf{u}}}&\enspace \sum_{i\in\mathcal{I}} \max(u_i,0)  \Delta T\ \\
\text{subject to:}\enspace&\eqref{eq:general_problem_discreteB}-\eqref{eq:general_problem_discreteC}, 
\eqref{eq:general_problem_discreteF}-\eqref{eq:general_problem_discreteN}, \eqref{eq:combined}.
\end{aligned}
\end{equation}
Substituting \eqref{eq:combined} into the objective~\eqref{eq:general_problem_discreteA}, we eliminate the control variable $u_i$ resulting in 
\begin{equation}
\begin{aligned}
\displaystyle{\minimize_{\mathbf{a}}}&\enspace \sum_{i\in\mathcal{I}} \max(a_{i} +d_1 + d_2 v_i + d_3 v_i^2,0)  \Delta T\ \\
\text{subject to:}\enspace&\eqref{eq:general_problem_discreteB}-\eqref{eq:general_problem_discreteC},\eqref{eq:general_problem_discreteF}-\eqref{eq:general_problem_discreteN}, 
\end{aligned}
\end{equation}
where the main optimization variable is now $\mathbf{a}$.

Since \( d_3 > 0 \) always holds, the function \( a_i + d_1 + d_2 v_i + d_3 v_i^2 \) is convex in the variables \( (a_i, v_i) \). Consequently, the overall problem is convex, as the objective function involves the pointwise maximum of two convex functions, and all remaining constraints are linear. To handle the $\max$ operator, we apply the standard epigraph reformulation (see~\cite{boyd}) by introducing auxiliary variables \( z_i \), yielding the following equivalent convex formulation referred to as ECO+ (Energy Consumption Optimization via Positive Control)
\begin{equation}
\label{eq:convex_form_eq}
\begin{aligned}
\displaystyle{\minimize_{\mathbf{a}}}&\enspace \sum_{i\in\mathcal{I}} z_i  \Delta T\ \\
\text{subject to:}\enspace&\eqref{eq:general_problem_discreteB}-\eqref{eq:general_problem_discreteC}, \eqref{eq:general_problem_discreteF}-\eqref{eq:general_problem_discreteN}, \\
& a_i + d_1 + d_2 v_i + d_3 v_i^2 \leq z_i,  \forall i \in \mathcal{I}^+, \\
& 0 \leq z_i,  \forall i \in \mathcal{I}^+,
\end{aligned}
\end{equation}
ECO+ is a Quadratically Constrained Program (QCP), which can be efficiently solved using modern convex optimization solvers. Since the control input \( \mathbf{u} \) was eliminated during reformulation, it can be recovered post-optimization using equation~\eqref{eq:combined}.

\begin{figure}[t]
	\centering
    \includegraphics[width=1\columnwidth]{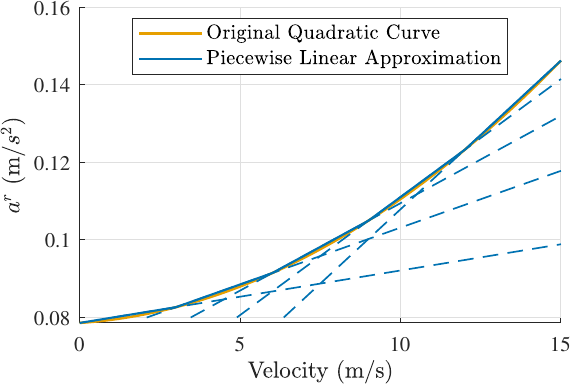}
	\caption{PWA affine approximation for the Quadratic constraint using CPEM model parameters.}
	\label{fig:Friction_approximation}
\end{figure}

\subsection{Quadratic Constraints Linearization}

The presence of quadratic constraints can substantially increase computational complexity. To mitigate this, we approximate the nonlinear term using piecewise affine (PWA) functions, offering a tractable alternative with minimal loss in accuracy. Specifically, we approximate the resistive force term, \( a^r_i \), which encapsulates the effects of aerodynamic drag, rolling resistance, and road slope, as
\begin{equation}
\label{eq:affine_approx}
a^r_i \leq y_k(v_i), \quad \forall k \in \{1, \dots, K\}, ~\forall i \in \mathcal{I}^+
\end{equation}
where \( y_k(v_i) = b^1_k v_i + b^2_k \) defines the \( k \)-th affine function, and \( K \) is the number of segments used in the approximation over the domain of \( v_i \in [0, v^{\max}] \).

To prevent underestimation, which would necessitate additional control input post-optimization, we exploit the convexity of the original constraint and construct affine segments by connecting adjacent evaluation points of \( a^r_i \) as 
\begin{equation}
\label{eq:affine_approx_coeff}
\begin{aligned}
b^1_k &= \frac{a^r_i(k\delta v) - a^r_i((k-1)\delta v)}{\delta v}, \\
b^2_k &= a^r_i((k-1)\delta v) - b^1_k \cdot (k-1)\delta v.
\end{aligned}
\end{equation}
where \( \delta v = v^{\max}/K \). This approach guarantees a conservative and continuous over-approximation. Additionally, it ensures continuity between adjacent segments and provides a smooth approximation across the entire velocity domain. Although the approximation must be applied for each index \( i \in \mathcal{I} \), the coefficients \( b^1_k \) and \( b^2_k \) are computed only once, yielding significant computational savings during optimization.

Figure~\ref{fig:Friction_approximation} illustrates the PWA approximation of the quadratic resistive force constraint using five linear segments for CPEM model parameters. The visual comparison shows that the approximation closely matches the true quadratic expression, with minimal deviation across the velocity range.

\section{Simulation Results and Discussion}
\label{Sec:Results}

In this section, we present and analyze the simulation results. We begin by detailing the simulation setup, followed by an evaluation of the accuracy of the PWA approximation. The benchmark strategies and performance metrics used for comparison are then described. We assess the proposed method against these baselines in terms of solution quality and computational efficiency. Finally, we examine the impact of a leading vehicle on the AV’s trajectory and investigate the effects of strict passenger comfort constraints.

\subsection{Simulation Setup}
\label{subsec:simulation_setup}

All simulations were conducted on a flat road segment of length \( L = 100 \) $\mathrm{m}$, with initial velocity \( v^o = 8 \) $\mathrm{m/s}$ and terminal velocity \( v^d \in \{6, 8, 10\} \) $\mathrm{m/s}$. Unless otherwise stated, the following parameters were used: \( v^{\min} = 0 \) $\mathrm{m/s}$, \( v^{\max} = 15 \) $\mathrm{m/s}$, \( u^{\min} = -3.5 \) $\mathrm{m/s}$\(^2\), \( u^{\max} = 2.5 \) $\mathrm{m/s}$\(^2\), \( J^{\min} = -10 \) $\mathrm{m/s}$\(^3\), \( J^{\max} = 10 \) $\mathrm{m/s}$\(^3\), and a discretization step size \( \Delta T = 0.1 \) $\mathrm{s}$.

The CPEM model parameters, derived from a Nissan Leaf, are summarized as follows: \( c_1 = 0.0328 \), \( c_2 = 4.575 \), \( C_r = 1.75 \), \( \rho = 1.2256 \) $\mathrm{kg/m}$\(^3\), \( g = 9.8066 \) $\mathrm{m/s}$\(^2\), \( A = 2.3316 \) $\mathrm{m}$\(^2\), and \( C_D = 0.28 \) . The drivetrain efficiency parameters are \( \eta_{\text{d}} = 0.92 \), \( \eta_{\text{em}} = 0.91 \), and \( \eta_{\text{b}} = 0.9 \). The regenerative braking efficiency is modeled using an exponential fit to empirical data as $\eta_{\text{rb}} = 0$ for $a(t) \geq 0$, and $\eta_{\text{rb}} = \left[\exp\left(0.0411 / |a(t)|\right)\right]^{-1}$ for $a(t) < 0$ \cite{CPEM}.

The KMMK model parameters are based on the torque-speed-efficiency map of a Nissan March K11, with coefficients: \( c_0 = 0.1569 \), \( c_1 = 0.0245 \), \( c_2 = -7.415 \times 10^{-4} \), \( c_3 = 5.975 \times 10^{-5} \), \( c_4 = 0.07224 \), \( c_5 = 0.09681 \), and \( c_6 = 1.075 \times 10^{-3} \). Additional constants include vehicle mass \( M = 1200 \) $\mathrm{kg}$, air density \( \rho = 1.184 \) $\mathrm{kg/m}$\(^3\), drag coefficient \( C_D = 0.32 \), frontal area \( A = 2.5 \) $\mathrm{m}$\(^2\), and rolling resistance coefficient \( \mu = 0.015 \) \cite{kamal}.

Each simulation uses friction coefficients consistent with the model being evaluated—that is, CPEM coefficients for the CPEM model and KMMK coefficients for the KMMK model.

All simulations were executed in MATLAB R2024b using the Gurobi Optimizer v12 \cite{gurobi} on a personal computer equipped with an AMD Ryzen 1600 processor and 16 GB of DDR4 RAM.

\begin{figure}[t]
	\centering
    \includegraphics[width=1\columnwidth]{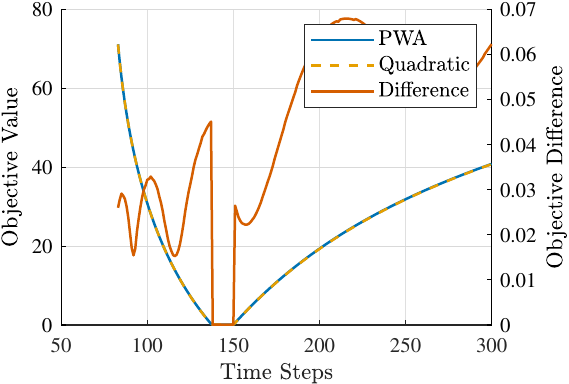}
	\caption{Comparison of objective values for varying time steps between PWA and Quadratic friction integrations.}
	\label{fig:Friction_comparison_cost}
\end{figure}
\begin{figure}[t]
	\centering
    \includegraphics[width=1\columnwidth]{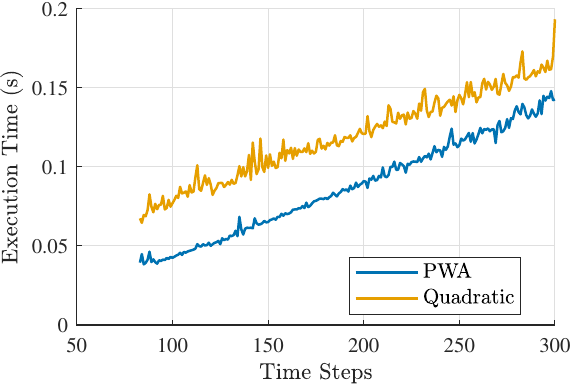}
	\caption{Execution time comparison for varying time steps with PWA and Quadratic constraint integrations.}
	\label{fig:Friction_comparison_execution_time}
\end{figure}

\begin{figure*}[t!]
	\centering
	\subcaptionbox{CPEM, $v^d=6$}[.3\linewidth][c]{%
		\includegraphics[width=1.0\linewidth]{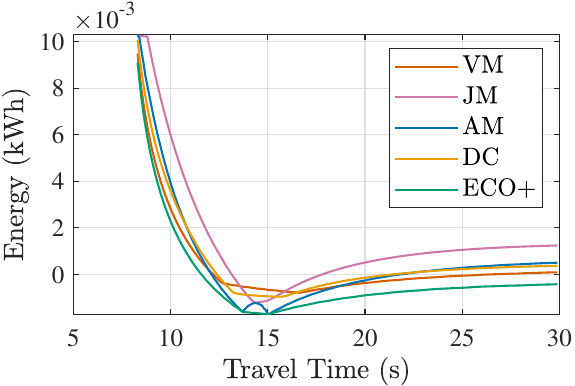}}\quad
	\subcaptionbox{CPEM, $v^d=8$}[.3\linewidth][c]{%
		\includegraphics[width=1.0\linewidth]{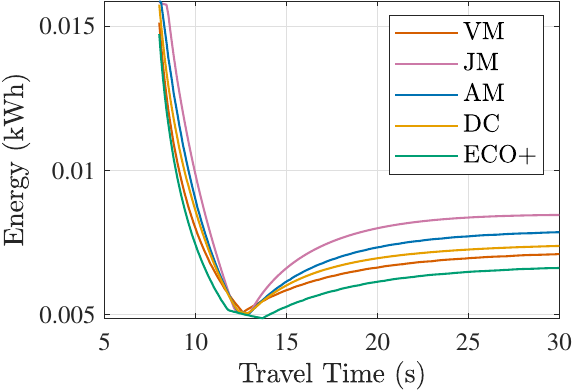}}\quad
	\subcaptionbox{CPEM, $v^d=10$}[.3\linewidth][c]{%
		\includegraphics[width=1.0\linewidth]{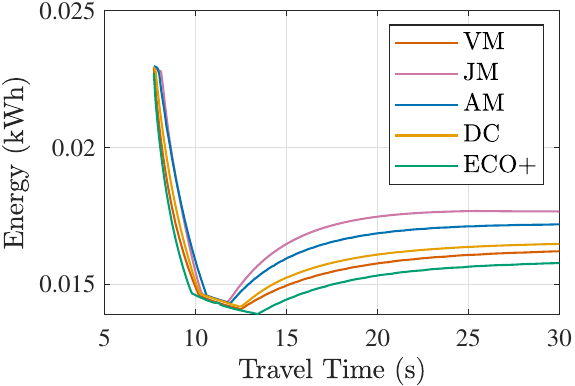}}		\\
	\subcaptionbox{KMMK, $v^d=6$}[.3\linewidth][c]{%
		\includegraphics[width=1.0\linewidth]{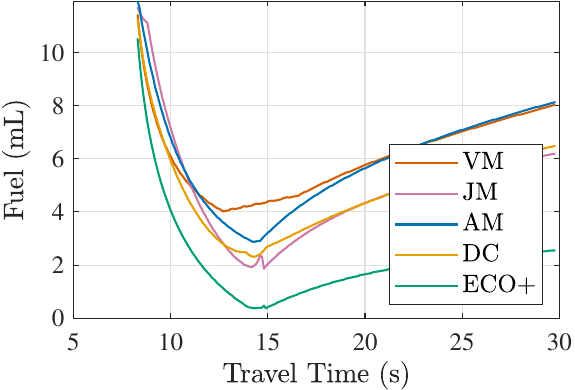}}\quad
	\subcaptionbox{KMMK, $v^d=8$}[.3\linewidth][c]{%
		\includegraphics[width=1.0\linewidth]{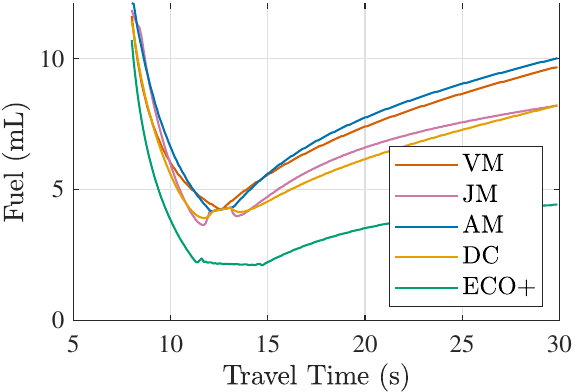}}\quad
	\subcaptionbox{KMMK, $v^d=10$}[.3\linewidth][c]{%
		\includegraphics[width=1.0\linewidth]{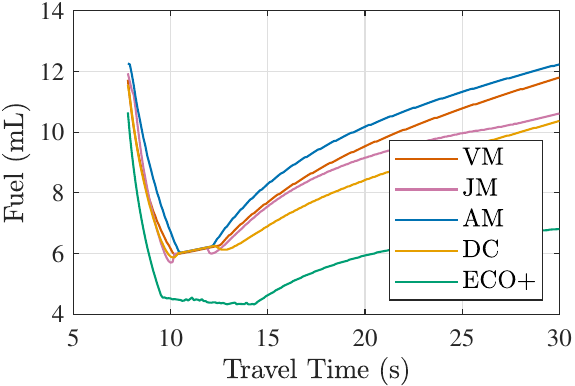}}			        
	\caption{Fuel and energy consumption as functions of travel time for the proposed ECO+ method and benchmark strategies (VM, JM, AM, DC).}
	\label{fig:Heuristic_comparison}
\end{figure*}

\subsection{Evaluation of the PWA Approximation}
\label{subSec:Results_pwa}

To quantify the accuracy of this approximation the PWA approximation, we consider a scenario with a desired terminal velocity of \( v^d = 6\,\text{m/s} \). Parameters from the CPEM model are employed in this experiment. The resulting objective values, computed across a range of travel times up to 30 $\mathrm{s}$, are shown in Figure~\ref{fig:Friction_comparison_cost}, which compares the performance of the exact quadratic constraint with its PWA approximation using five linear segments.

The difference in objective value between the two methods is negligible. On average, the PWA-based solution differs from the quadratic formulation by less than 0.4\%, indicating that the approximation retains high fidelity to the optimal behavior.

The computational benefits of the PWA formulation are demonstrated in Figure~\ref{fig:Friction_comparison_execution_time}, which reports execution times for both approaches. As expected, execution time increases with travel duration due to the expanding optimization horizon. Nonetheless, the PWA formulation consistently outperforms the quadratic version, achieving an average speed-up of 36.9\%.

In summary, the PWA approximation provides a substantial reduction in computation time with no significant loss in solution quality. Similar results have been observed when applying the KMMK model parameters and across various initial and terminal velocity combinations.

\subsection{Benchmark Strategies for  Consumption Optimization}
\label{subsec:baseline-methods}

This subsection presents several benchmark strategies for minimizing consumption. Each strategy is characterized by a distinct objective function, while the underlying optimization framework remains consistent across all cases. These strategies serve as reference points for evaluating the effectiveness of the proposed approach.

\begin{enumerate}
    \item \textbf{Velocity Minimization (VM):} This strategy minimizes the squared $\mathcal{L}^2$ norm of the velocity,
    \[
    \left\lVert v \right\rVert^2_2 = \int_{t^o}^{t^d} v^2(t)\, dt,
     \]
    which corresponds to minimizing the vehicle’s kinetic energy. Minimizing kinetic energy reduces the propulsion effort required, thereby improving energy efficiency. As such, this strategy is frequently employed in eco-driving frameworks to limit deviations from a reference velocity \cite{v_cubed,v_cubed_cas}.
    
    \item \textbf{Jerk Minimization (JM):} To enhance passenger comfort, this strategy minimizes the squared $\mathcal{L}^2$ norm of jerk,
    \[
    \left\lVert J \right\rVert^2_2 = \int_{t^o}^{t^d} J^2(t)\, dt.
     \]
    Minimizing jerk promotes smoother transitions in acceleration and deceleration, improving ride quality \cite{Zhang2019}.
    
    \item \textbf{Acceleration Minimization (AM):} This strategy minimizes the squared $\mathcal{L}^2$ norm of acceleration,
    \[
    \left\lVert a \right\rVert^2_2 = \int_{t^o}^{t^d} a^2(t)\, dt,
     \]
    a formulation shown to reduce transient engine operation and thereby improve fuel efficiency~\cite{Malikopoulos2013}. It is widely used as a surrogate for consumption \cite{Malikopoulos2018,xu2022general,ahadjistimoITV,Cas_barrier}.
    
    \item \textbf{DC programming (DC):} This approach extends~\cite{fuel_convex_concave} by approximating the nonlinear consumption model with a quadratic surrogate,
    \[
    c_1 a_i^2 + c_2 v_i^2 + c_3 a_i v_i + c_4 a_i + c_5 v_i + c_6,
    \]
    where the coefficients are fitted to the original model. The resulting problem is solved using Difference-of-Convex (DC) programming.
    
    \item \textbf{Nonlinear Solver (NLS):} This strategy directly optimizes the nonlinear consumption model using MATLAB’s \textit{fmincon} solver. Multiple initial conditions are tested to reduce sensitivity to local minima.\footnote{The options for fmincon were set to MaxFunctionEvaluations = $10^7$ and MaxIterations = $10^5$.}
\end{enumerate}

Alternative strategies, including the $\mathcal{L}^2$ norm of the control input, the corresponding $\mathcal{L}^1$ norms, and weighted linear combinations such as $\left( v(t) + a(t)^2 \right)dt$, were explored but consistently exhibited inferior performance compared to the aforementioned formulations. Consequently, they were excluded from further analysis.

To evaluate solution quality, we report the relative difference in consumption between methods, defined as
\begin{equation}
\label{eq:relative_difference}
\text{Relative Difference} = \frac{|x - y|}{\max(|x|, |y|)} \times 100\%,
\end{equation}
where \( x \) and \( y \) denote the energy costs associated with two different strategies.

\begin{figure}[t]
	\centering
	 \begin{subfigure}[b]{\columnwidth}
	 	\centering
	 	\includegraphics[width=1\linewidth]{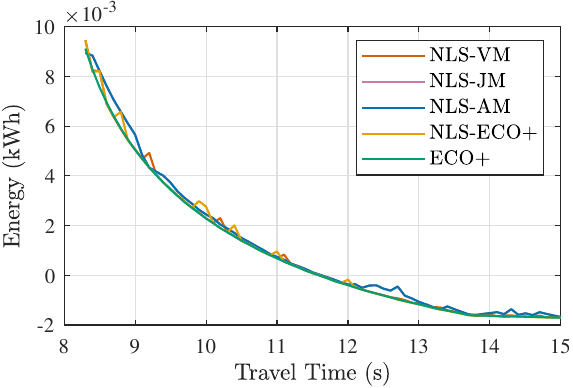}
         	\caption{CPEM model}
	 	\label{fig:CPEM_Advanced_comp}
	 \end{subfigure}%
	
	 \begin{subfigure}[b]{\columnwidth}
		\centering
		\includegraphics[width=1\linewidth]{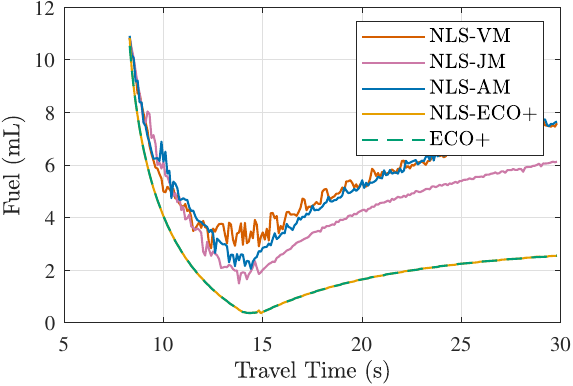}
				 \caption{KMMK Model}
		 \label{fig:KMMK_Advanced_comp}
	 \end{subfigure}
	\caption{Fuel and energy consumption as functions of travel time for the proposed ECO+ method and the NLS initialized with VM, JM, AM, and ECO+ solutions.}
	\label{fig:Advanced_comp}
\end{figure}

\textbf{Remark:} ECO+ is implemented with the PWA approximation for the remainder of this paper.

\subsection{Objective Performance and Execution Time}
\label{subSec:Results_objective_execution}

Several benchmark strategies were introduced in Section~\ref{subsec:baseline-methods}. To evaluate the performance of the proposed ECO+ approach, we consider three scenarios, each involving a single CAV exiting with final speeds of \(v^d  \) equals to 6, 8, and 10 $\mathrm{m/s}$. For each scenario, the travel time \(t^m\) is varied in increments of \(\delta t = 0.1\) s to construct trade-off curves between consumption and travel time.

First, ECO+ is compared against VM, JM, AM, and DC, all implemented using standard optimization techniques. Figures~\ref{fig:Heuristic_comparison}~(a)--(c) and (d)--(f) show the minimum consumption attained by each strategy for the CPEM and KMMK models, respectively. In every case, consumption is high for small travel times, \(t^m\). However, as \(t^m\) increases, consumption decreases until a minimum is reached. Beyond this point, further increases in \(t^m\) result in a monotonic rise in fuel consumption. For \(v^d = 6\)~$\mathrm{m/s}$ in the CPEM model (Figure~\ref{fig:Heuristic_comparison}(a)), the net energy occasionally becomes negative, indicating energy recuperation via regenerative braking during extended deceleration.

\begin{figure}[t]
	\centering

		\includegraphics[width=1\linewidth]{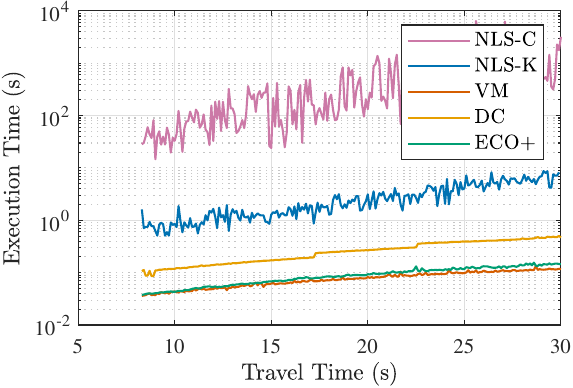}

	\caption{Execution time as a function of travel time for VM, DC, ECO+, and NLS methods. NLS-C and NLS-K correspond to the direct optimization of the CPEM and KMMK models, respectively, highlighting differences in problem complexity.}
	\label{fig:Advanced_execution}
\end{figure}

VM consistently achieves lower energy consumption than JM but yields higher fuel consumption, suggesting that different strategies better suit electric versus conventional powertrains. Across all travel times and scenarios, ECO+ attains the lowest consumption for both models.

Using the relative difference metric in~\eqref{eq:relative_difference}, ECO+ outperforms VM by 69.6\%, 7.3\%, and 2.6\% for \(v^d = 6\), 8, and 10~$\mathrm{m/s}$ in the CPEM model, respectively. In the KMMK model, ECO+ surpasses DC by reducing consumption by 50.6\%, 40.4\%, and 29.8\% for \(v^d = 6\), 8, and 10~$\mathrm{m/s}$, respectively. The performance gap narrows at higher speeds and shorter travel times as the feasible set is small and widens for longer durations as the feasible set expands.

Next, ECO+ is benchmarked against the nonlinear solver (NLS) initialized with VM, JM, AM, and ECO+ solutions, denoted as NLS-VM, NLS-JM, NLS-AM, and NLS-ECO+. DC was excluded due to its high computational cost. Figure~\ref{fig:Advanced_comp} presents results for \(v^d = 6\)~$\mathrm{m/s}$. For the CPEM model, ECO+ remains superior, as NLS fails to find a better solution and frequently exhausts the evaluation limit. However, all NLS variants converge near the ECO+ solution, suggesting that ECO+ approximates a local, or potentially global, optimum. 
For the KMMK model, NLS-ECO+ demonstrates a marginal 0.27\% improvement over ECO+ by more precisely satisfying the quadratic constraint~\eqref{eq:general_problem_discreteE}. In contrast, when initialized with alternative solutions, the nonlinear solver produces significantly poorer results. This minor improvement highlights that ECO+ offers a high-quality approximation.

Figure~\ref{fig:Advanced_execution} reports execution times, including the construction of constraint matrices. As the travel time increases, the execution time increases for all approaches, since the horizon and hence the number of variables increases. ECO+ requires solving a single linear program (LP), with an average runtime of 0.09~$\mathrm{s}$. VM, AM, and JM, formulated as quadratic programs, are on average 17.7\% faster because they omit the nonlinear resistive-force constraint~\eqref{eq:general_problem_discreteE}. At shorter time horizons, ECO+ and VM runtimes are similar, but for longer time horizons, the cost of matrix construction accumulates. DC is roughly three times slower than ECO+, averaging 0.27~$\mathrm{s}$. The NLS requires 2.7~$\mathrm{s}$ for the KMMK model and frequently fails to converge within limits for the CPEM model, averaging 746~$\mathrm{s}$.

In summary, ECO+ consistently delivers the best solution quality across all scenarios, while maintaining execution times close to the fastest approach. ECO+ achieves superior performance by effectively modeling resistive forces, enabling the vehicle to glide without active control over extended periods, akin to the PnG strategy.

\subsection{Leading vehicle}
\label{subSec:Results_leading}

\begin{figure}[t]
	\centering
    \includegraphics[width=1\columnwidth]{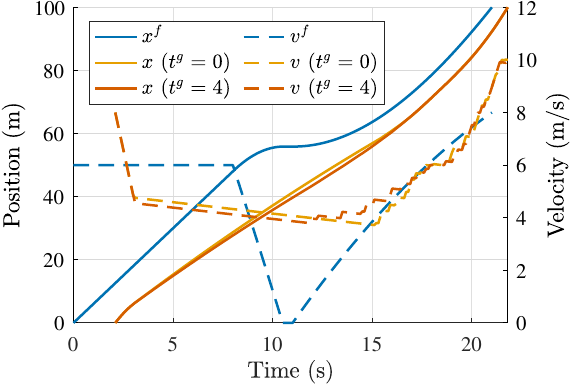}
	\caption{Trajectory and velocity profiles of the preceding vehicle and the ego vehicle with and without the time gap constraint.}
	\label{fig:Position_Follow}
\end{figure}

	\begin{figure}[t]
	\centering
	\begin{subfigure}[b]{\columnwidth}
		\centering
		\includegraphics[width=1\linewidth]{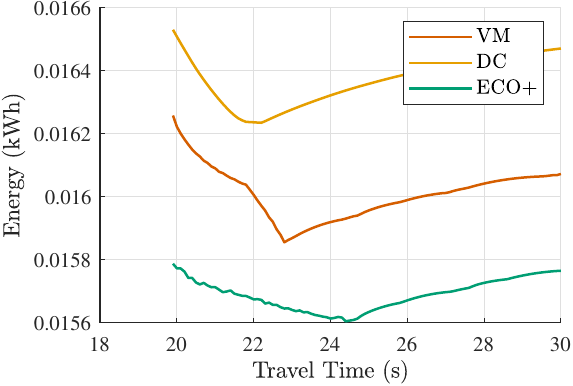}
		\caption{CPEM model}
		\label{fig:Tradeoff_follow_CPEM}
	\end{subfigure}%
	
	\begin{subfigure}[b]{\columnwidth}
		\centering
		\includegraphics[width=1\linewidth]{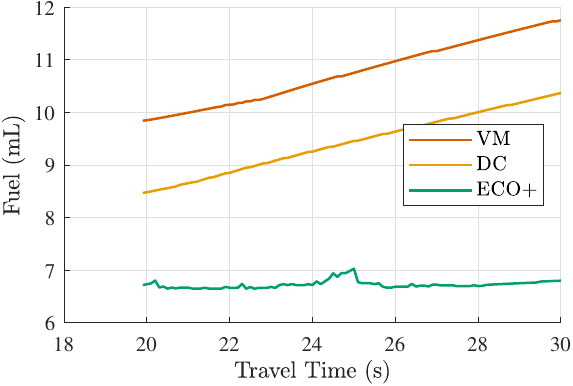}
				\caption{KMMK Model}
		\label{fig:Tradeoff_follow_KMMK}
	\end{subfigure}
	\caption{Fuel and energy consumption of the ego vehicle as a function of travel time for the VM, DC, and ECO+ methods, considering the presence of a preceding vehicle.}
	\label{fig:Tradeoff_follow}
\end{figure}

In the previous results, we considered a single AV scenario. Here, we examine the impact of having multiple vehicles on consumption. Specifically, we introduce a leading vehicle following a predefined velocity profile, $v^f$, as shown in Figure \ref{fig:Position_Follow}. The leading vehicle enters the zone with a speed of 6 \(\mathrm{m/s}\), suddenly stops at 10.5 $\mathrm{s}$, remains stationary for 0.5 $\mathrm{s}$, and then accelerates, minimizing $\mathcal{L}^2$ norm of acceleration, to an exit velocity of 8 \(\mathrm{m/s}\) at 21 $\mathrm{s}$. This velocity profile was chosen to significantly impact the trajectory of the following AV, as the sudden stop and relatively slow acceleration force the following vehicle to adjust its movement to maintain a safe distance.

The ego AV enters the zone after a delay of 2 $\mathrm{s}$ with an initial velocity of \(v^o = 8\, \mathrm{m/s}\) and exits with a final velocity of \(v^d = 10\, \mathrm{m/s}\). As in the previous simulations, the travel time is increased up to 30 $\mathrm{s}$ to construct a tradeoff curve between travel time and energy consumption. It is important to note that the inclusion of a leading vehicle introduces additional constraints, such as maintaining a safety distance and a time gap, as defined in the problem formulation. To illustrate the effect of the time gap constraint, Figure~\ref{fig:Position_Follow} compares two scenarios: \(t^g = 0\) and \(t^g = 4\) $\mathrm{s}$, with the ego AV’s travel time set to \(t^m = 19.7\) $\mathrm{s}$. As the leading vehicle comes to a stop, the following vehicle adjusts by reducing its velocity earlier to maintain the required safety distance, thereby altering its acceleration profile.

The corresponding tradeoff curves are shown in Figure~\ref{fig:Tradeoff_follow} for the VM, DC, and ECO+ solution approaches with \(t^g = 4\) $\mathrm{s}$. The results are consistent with the trends observed in the previous section, confirming that the relative performance ranking among the methods is preserved. Compared to ECO+, the VM and DC approaches result in, on average, 2\% and 4.1\% higher consumption under the CPEM model, and 37.54\% and 28.77\% higher consumption under the KMMK model, respectively.

It can also be observed that the difference in performance between the methods is considerably less compared to the single CAV scenario. This is because the feasible acceleration profiles of the following vehicle are constrained by the trajectory of the leading vehicle. Thus, it becomes apparent that the most efficient scenario would involve a control zone exclusively for AVs, as this would allow for uniform acceleration profiles and optimized for consumption.

\subsection{Impact of Tight Comfort Constraints}
\label{subSec:Results_comfort_tradeoff}

\begin{figure}[t]
	\centering
    \includegraphics[width=1\columnwidth]{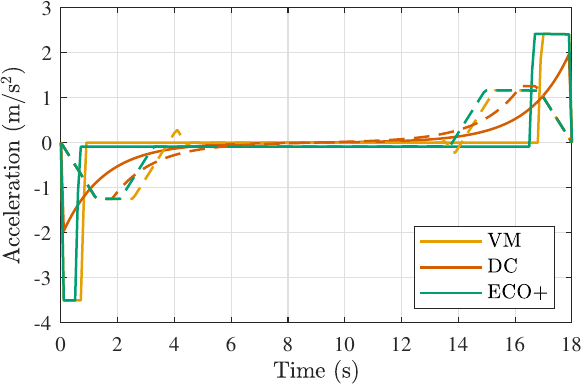}
	\caption{Acceleration profiles of the VM, DC, and ECO+ methods, before and after applying comfort constraints.}
	\label{fig:Profiles}
\end{figure}
	\begin{figure}[t]
	\centering
	\begin{subfigure}[b]{\columnwidth}
		\centering
		\includegraphics[width=1\linewidth]{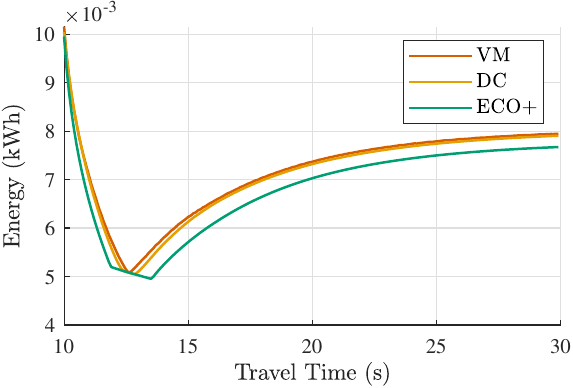}
				\caption{CPEM}
		\label{fig:Tradeoff_jerk_CPEM}
	\end{subfigure}	%
	\begin{subfigure}[b]{\columnwidth}
		\centering
		\includegraphics[width=1\linewidth]{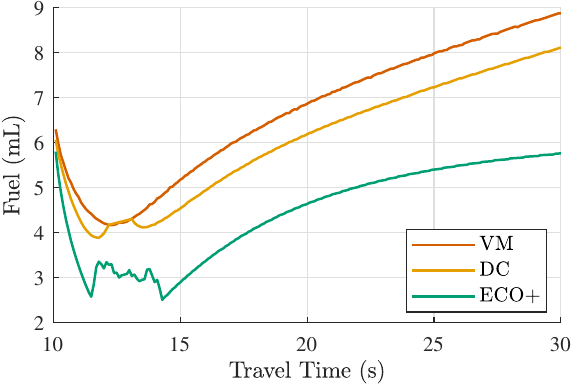}
				\caption{KMMK}
		\label{fig:Tradeoff_jerk_KMMK}
	\end{subfigure}
	\caption{Fuel and energy consumption as a function of travel time for the VM, DC, and ECO+ methods, under strict comfort constraints.}
	\label{fig:Tradeoff_jerk}
\end{figure}

Each method produces unique acceleration profiles based on the optimization objective. Example acceleration profiles for a travel time of 18 $\mathrm{s}$ and a merging speed of 8 $\mathrm{m/s}$, resulting from the VM, DC, and ECO+ solution approaches, are shown in Figure \ref{fig:Profiles}. Initially, all approaches involve deceleration. The VM approach applies maximum deceleration for approximately 0.7 $\mathrm{s}$, followed by cruising at a constant speed. At 17 $\mathrm{s}$, it initiates maximum acceleration for 1 second to reach the merging speed. The ECO+ approach also begins with maximum deceleration, but over a shorter duration (0.5 $\mathrm{s}$), after which the vehicle enters a gliding phase. At 16.5 $\mathrm{s}$, it starts accelerating at the maximum allowable rate to achieve the required speed. In contrast, the DC approach gradually reduces deceleration and transitions to acceleration earlier, beginning around 9 $\mathrm{s}$, with the acceleration rate progressively increasing, particularly after 15 $\mathrm{s}$.

Rapid changes in acceleration, i.e., jerk, and high levels of acceleration or deceleration are associated with passenger discomfort \cite{comfort-hoberock1976survey}. The ISO 15622:2018 standard for adaptive cruise control specifies that the average acceleration over one second should be limited to 2 $\mathrm{m/s}^2$, while deceleration should not exceed -3.5 $\mathrm{m/s}^2$ \cite{ISO15622:2018}. Additionally, it mentions a limit on negative jerk of -2.5 $\mathrm{m/s}^3$. However, experimental studies, such as \cite{Comfort1}, indicate that acceptable levels for passenger comfort are even lower.

For this reason, we introduce stricter limits on both jerk and acceleration: $J^\min = -1  ~\mathrm{m/s}^3$
 , $J^\max = 1 ~\mathrm{m/s}^3$
 , and $a^\min = -1.25  ~\mathrm{m/s}^2$
 , $a^\max = 1.25  ~\mathrm{m/s}^2$, based on findings from \cite{Comfort1}. The changes in the acceleration profiles obtained for solving Problem \eqref{eq:general_problem_discrete} under these stricter limits are shown in Figure \ref{fig:Profiles} with dashed lines. Although the overall profile patterns remain similar, all methods now exhibit smoother deceleration and acceleration, which is consistent with the anticipated effects of the imposed jerk constraint.

We construct the tradeoff between travel time and consumption, which is displayed in Figure \ref{fig:Tradeoff_jerk} for both models. Although the relative performance of the methods shifts under tighter comfort constraints, the ECO+ approach consistently outperforms the alternatives. However, the performance gap is notably reduced. On average, the VM and DC solution approaches result in higher consumption: 4.7\% and 4\% for the CPEM model, and 33\% and 26.1\% for the KMMK model, compared to the ECO+ solution approach. Notably, before applying the stricter limits, the performance gap between the DC and ECO+ approaches was 10.9\% and 40.4\% for the CPEM and KMMK models, respectively.

The reduction in the performance gap can primarily be attributed to the stricter limits on jerk and acceleration, which significantly reduce the feasible region of the solution space. While the limits used in these simulations are stricter than those specified by ISO standards, they still demonstrate the effectiveness of the ECO+ strategy in minimizing consumption while adhering to strict comfort constraints.

\section{Conclusion}
\label{sec:conclusion}

This work introduced the PCI novel surrogate consumption model that considers positive control input, and integrated it into ECO+ (Energy Consumption Optimization via Positive Control), a convex, time-based optimization framework for autonomous vehicles approaching intersections. The proposed approach ensures safety and passenger comfort while optimizing energy use. To enable scalable real-time implementation, we approximate quadratic constraints using piecewise affine segments, thereby reformulating the problem as a linear program without compromising solution quality.

ECO+ was validated using two distinct empirical models, one for fuel consumption and one for electric energy consumption, and benchmarked against various optimization strategies, including a state-of-the-art nonlinear solver. Extensive simulation results demonstrated that ECO+ consistently outperforms alternative methods in reducing both fuel and energy consumption, even under strict comfort constraints and in scenarios involving a leading vehicle. 

While the simulation results are promising, the absence of experimental validation remains a limitation. Real-world testing presents logistical challenges, including the need for precise instrumentation, controlled environments, and coordinated multi-vehicle deployments. Future work should focus on experimental evaluation to assess the practical feasibility and robustness of ECO+ under realistic conditions, including sensor noise and environmental disturbances.

Overall, ECO+ offers a practical, scalable, and computationally efficient solution for improving the energy efficiency and sustainability of autonomous urban mobility systems.

	\bibliographystyle{IEEEtran}
	\bibliography{bibliography}
	
		\begin{IEEEbiography}[{\includegraphics[width=1in,height=1.25in,clip,keepaspectratio]{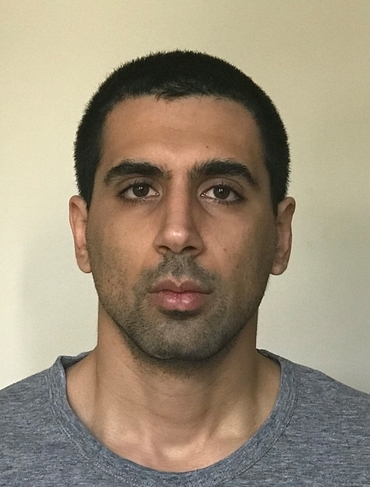}}]{Andreas Hadjigeorgiou} received his BEng in Mechanical Engineering from the University of Bath and his BSc in Applied Mathematics from the University of Cyprus in 2012 and 2016 respectively. Currently, he is a Ph.D. candidate at the Department of Electrical and Computer Engineering of the University of Cyprus.  His research interests  include  optimization   techniques, Intelligent Transportation Systems and Connected and Autonomous vehicles.
		\end{IEEEbiography}
        \begin{IEEEbiography}[{\includegraphics[width=1in,height=1.25in,clip,keepaspectratio]{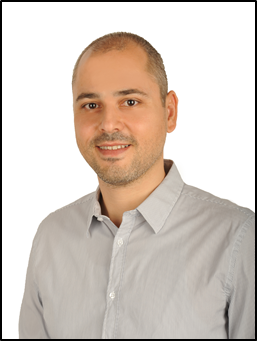}}]
        {Stelios Timotheou} received the Dipl.-Ing. degree in electrical and computer engineering from the National Technical University of Athens and the M.Sc. degree in communications and signal processing and the Ph.D. degree in intelligent systems and networks from the Department of Electrical and Electronic Engineering, Imperial College London, in 2010. In previous appointments, he was a Research Associate with KIOS, Department of Electrical and Computer Engineering, University of Cyprus, as a Visiting Lecturer, and a Post-Doctoral Researcher with the Computer Laboratory, University of Cambridge. He is currently an Associate Professor with the Department of Electrical and Computer Engineering and a Faculty Member with the KIOS Research and Innovation Center of Excellence, University of Cyprus. His research focuses on monitoring, control, and optimization of critical infrastructure systems, with emphasis on intelligent transportation systems and communication systems. He is a Senior Editor of IEEE Transactions on Intelligent Transportation Systems and an Associate Editor of IEEE Transactions on Intelligent Vehicles.

        \end{IEEEbiography}\vfill

\end{document}